\documentclass[A4]{article}
\usepackage{color}
\usepackage{latexsym}
\usepackage{graphicx}
\usepackage{amssymb,amsmath}
\usepackage{authblk}

\usepackage[latin1]{inputenc}

\newtheorem{thm}{Theorem}
\newtheorem{prop}[thm]{Proposition}

\makeatletter
\let\@fnsymbol\@arabic
\makeatother

\newcommand{\R}{\mathbb R}
\newcommand{\erre}{\mathbb R}

\newcommand{\enne}{\mathbb N}
\newcommand{\errepz}{\mbox{$\erre^{+}_0$}}

\newcommand{\Rp}{\mathbb R_+}
\newcommand{\Ks}{\mbox{$\cal K$}}

\newcommand{\freccia}{\mbox{$\rightarrow$}}

\def\epsilon{\varepsilon}

\def\qed{\hskip 1mm\boxit{}\hskip 1mm}
\def\boxit#1{\vbox{\hrule\hbox{\vrule\kern.75truemm
\vbox{\kern.75truemm#1\kern1truemm}\kern1truemm\vrule}\hrule}}

\title{ 
Convergence and rate of approximation in $BV^{\varphi}(\R^N_+)$ for a class of Mellin integral operators
}
\author{Laura Angeloni \footnote{Laura Angeloni - email: laura.angeloni@unipg.it, Phone: +39 075 585 5036, Fax: +39 075 585 5024}
\\
Gianluca Vinti \footnote{Gianluca Vinti (corresponding author) - email: gianluca.vinti@unipg.it, Phone: +39 075 585 5025, Fax: +39 075 585 5024} \\ \vskip0.5cm  
\small \it Dipartimento di Matematica e Informatica \\ \it Universit\`{a} degli Studi di Perugia,\\  \small \it Via Vanvitelli 1, 06123 Perugia (Italy)}

\date{}

\begin{document}

\maketitle

\begin{abstract}
In this paper we study convergence results and rate of approximation for a family of linear integral operators of Mellin type in the frame of $BV^{\varphi}(\R^N_+)$. Here $BV^{\varphi}(\R^N_+)$ denotes the space of functions with bounded $\varphi-$variation on $\R^N_+$, defined by means of a concept of multidimensional $\varphi-$variation in the sense of Tonelli.  
\end{abstract}


\vskip0.2cm


\section{Introduction}\label{sec-intro}

The importance of Mellin operators in approximation theory is well-known: they are widely studied (see, e.g., \cite{MAM,BJ1}) and they have important applications in several fields. For example, we recall that Mellin analysis has deep connections with Signal Processing, in particular with the so-called Exponential Sampling (see \cite{BJ2}).

In this paper we study approximation properties for a family of linear integral operators of Mellin type of the form 
$$
(T_w f)({\tt s}) =\int_{\R^N_+} K_w({\tt t}) f({\tt st}) \langle{\tt t}\rangle^{-1} \,d {\tt t}, \ \ {\tt s}\in\R^N_+,\ w>0, \eqno {\rm(I)}
$$
with respect to the multidimensional $\varphi-$variation in the sense of Tonelli introduced in \cite{AVI7}. Here $\{K_w\}_{w>0}$ is a family of approximate identities (see Section \ref{sec-not}), $\langle{\tt t}\rangle:=\prod_{i=1}^N t_i$ and ${\tt s t}:=(s_1 t_1,\dots,s_N t_N)$, ${\tt s,t}\in\R^N_+$.

The class of the above operators (I) contains, as particular cases, several families of well-known integral operators (see Section \ref{sec-order}): among them, for example, the moment-type or average operators, the Mellin Picard operators and others. 

Due to the homothetic structure of our operators, it seems that the most natural way to frame the theory is to work with the Haar measure in $\R^N_+$, i.e., $\mu(A):=\int_A \langle{\tt t}\rangle^{-1}\,d {\tt t}$, where $A$ is a Borel subset of $\R^N_+$. Results about homothetic-type operators in various settings can be found, for example, in \cite{BV1,MAM,V,BSV,SCIAVI1,BMV,BM1,BM2,A2,AVI7,AVI8,AVI9}, while for similar results about classical convolution operators see, e.g., \cite{BUNE,SZ,MAVI,BABUSTVI,BV,AVI3,AVI4,AVI5,A1,A3}.

The main results are presented in Sections \ref{sec-est} and \ref{sec-order}. We first study the problem of the convergence in $\varphi-$variation: in particular, after some estimates for our integral operators, we prove that, if $f\in AC^{\varphi}(\R^N_+)$ (the space of $\varphi-$absolutely continuous functions), there exists a constant $\mu>0$ such that 
$$
\lim_{w\to +\infty} V^{\varphi}[\mu(T_w f-f)]=0.\eqno {\rm(II)}
$$
Then we face the problem of the rate of approximation and we prove that, if $f$ belongs to a Lipschitz class $V^{\varphi}Lip_{N}(\alpha)$, $\alpha>0$, under suitable assumptions on the kernels $\{K_w\}_{w>0}$ (see Section \ref{sec-order}), there exists a constant $\lambda>0$ such that 
$$
V^{\varphi}[\lambda(T_w f-f)]=O(w^{-\alpha}),
$$
for sufficiently large $w>0$.

An important step in order to achieve (II) is to prove the convergence for the $\varphi-$modulus of smoothness in the present setting; this problem was solved in \cite{AVI7}. This result extends to the multidimensional case an analogous one for the (one-dimensional) Musielak-Orlicz $\varphi-$variation (\cite{SCIAVI1}). In the case of the classical variation (see, e.g., \cite{BABUSTVI} for translation operators) such result is an easy consequence of the integral representation of the variation for absolutely continuous functions; on the contrary, in the case of the $\varphi-$variation, due to the lack of an integral representation, it requires a more delicate direct construction.

\section{Notations}\label{sec-not}

We will study approximation results in $BV^{\varphi}(\R^N_+)$, namely the space of functions $f:\R^N_+\rightarrow \R$ of bounded $\varphi-$variation introduced in \cite{AVI7}.  Such a concept of multidimensional $\varphi-$variation on $\R^N_+$ has the purpose to provide a $\varphi-$variation in the sense of Musielak-Orlicz (\cite{MO}) in the multidimensional frame, following the Tonelli approach (\cite{TO}), generalized in dimension $N\ge 2$ by T. Rad\'o (\cite{RA}) and C. Vinti (\cite{VI}).  Here we endow $\R^N_+$ with the Haar measure $\mu(A)=\int_A \langle{\tt t}\rangle^{-1} \,d{\tt t}$, where $A$ is a Borel subset of $\R^N_+$, $\langle {\tt t}\rangle:=\prod_{i=1}^N t_i$, ${\tt t}=(t_1,\dots,t_N)\in\R^N_+$, which seems to be the natural setting working with homothetic operators.
We recall that, under some properties of approximate continuity, the multidimensional version in the sense of Tonelli of the classical variation is equivalent to the distributional variation (see, e.g., \cite{CE,DE,GIU}).

We denote by $\Phi$ the class of all the functions $\varphi:\R^+_0 \freccia \R^+_0$ such that
\begin{enumerate}
\item $\varphi$ is convex and $\varphi(u)=0$ if and only if $u=0$;
\item $u^{-1} \varphi(u) \freccia 0$ as $u\freccia 0^+.$
\end{enumerate}
From now on we will assume that $\varphi\in\Phi$.

We now recall some notations of the multidimensional setting in which we work (see, e.g., \cite{BABUSTVI}). 
For $f:\R_+^N\freccia\R$ and ${\tt
x}=(x_1,\dots,x_N)\in\Rp^N$, $N\in\enne$, if we want to focus the attention on 
the $j-$th coordinate, $j=1,\dots,N$, we will write
$$
{\tt x}'_j= (x_1,\dots, x_{j-1}, x_{j+1},\dots,x_N)\in\Rp^{N-1},\ \ 
{\tt x}=({\tt x}'_j,x_j),\ \ f({\tt x})=f({\tt x}'_j,x_j).
$$
Given $I=\prod_{i=1}^N
[a_i,b_i]\subset\R^N_+$, by $I'_j:=[{\tt a}'_j,{\tt b}'_j]$ we will denote the
\hbox{$(N-1)$-} dimensional interval obtained deleting by $I$ the
$j-$th coordinate, so that
$$
I=[{\tt a}'_j,{\tt b}'_j]\times[a_j,b_j].
$$

In order to define the multidimensional $\varphi-$variation, we first recall that the $\varphi-$variation of a function $g:[a,b]\freccia\erre$
is defined as
$$
V^{\varphi}_{[a,b]}[g] := \sup_{D} \sum_{i=1}^n
\varphi(|g(s_i)-g(s_{i-1})|),
$$
where $D=\{s_0=a,s_1,\dots,s_n=b\}$ is a partition of $[a,b]$  (\cite{MO,MU}), and $g$ is said to be of bounded $\varphi-$variation ($g\in BV^{\varphi}([a,b])$) if $V^{\varphi}_{[a,b]}[\lambda g]<+\infty$, for some $\lambda>0$.
The $\varphi-$variation was introduced by L.C. Young (\cite{YO2}) as a generalization of the concept of $p-$variation, $p\ge 1$ (\cite{YO1,LY}), which extends Wiener's quadratic variation (\cite{WI}). However the main developments of this concept are due to J. Musielak and W. Orlicz and their school: we refer to \cite{MO} for the main properties of the (one-dimensional) $\varphi-$variation. For results concerning the $\varphi-$variation, the reader can see, e.g., \cite{MO,HER,MU,MAOR,MATOR,RAMA,AD, CHGA, SCIAVI2}. 

Now we consider the Musielak-Orlicz $\varphi-$variation of the $j-$th section of $f$, i.e., $V^{\varphi}_{[a_j,b_j]}[f({\tt x}'_j,\cdot)]$, for ${\tt x}'_j\in I'_j$, and then the $(N-1)$-di\allowbreak mensional integrals 
$$
\Phi^{\varphi}_j(f,I):= \int_{{\tt a}'_j}^{{\tt b}'_j} V^{\varphi}_{[a_j,b_j]}[f({\tt x}'_j,\cdot)]
{d{\tt x}'_j \over \langle {\tt x}'_j\rangle},
$$
where $\langle {\tt x}'_j\rangle:=\prod_{i=1, i\neq j}^N x_i$.

We now denote by
$$
\Phi^{\varphi}(f,I):= \left\{\sum_{j=1}^N [\Phi^{\varphi}_j(f,I)]^2\right\}^{1\over 2},
$$
the euclidean norm of $(\Phi^{\varphi}_1(f,I),\dots, \Phi^{\varphi}_N(f,I))$,
where we put $\Phi^{\varphi}(f,I)=+\infty$ if $\Phi^{\varphi}_j(f,I)=+\infty$ for some
$j=1,\dots,N$. Then the multidimensional $\varphi-$variation of $f$ on an interval $I\subset \Rp^N$ is defined as
$$
V_I^{\varphi}[f] := \sup \sum_{i=1}^m \Phi^{\varphi}(f,J_{i}),
$$
where the supremum is taken over all the finite families of
$N-$dimensional intervals $\{J_{1},\dots,J_{m}\}$ which form
partitions of $I$.\\
Finally by
$$
V^{\varphi}[f]:= \sup_{I\subset \Rp^N} V_I^{\varphi}[f],
$$
where the supremum is taken over all the intervals $I\subset \Rp^N$,
we will denote the $\varphi-$variation of $f$ over the whole space $\Rp^N$.

We will say that a function $f$ is {\it of bounded $\varphi-$variation} on $\R_+^N$ if there exists a constant $\lambda>0$ such that $V^{\varphi}[\lambda f]<+\infty$ and
$BV^{\varphi}(\R_+^N)$ will denote the space of functions of bounded $\varphi-$variation on $\R^N_+$, namely
$$
BV^{\varphi}(\R_+^N):=\{f\in {\mathcal M} :\ \exists\lambda>0\ \hbox{s.t.}\
V^{\varphi}[\lambda f]<+\infty\},
$$
where ${\mathcal M}$ is the space of all the measurable functions $f:\R^N_+\freccia \R$.
For the main properties of the multidimensional $\varphi-$variation, see \cite{AVI7}.

Finally by $AC^{\varphi}_{loc}(\R_+^N)$ we will denote the space of functions $f:\R_+^N\freccia \R$ which are {\it locally $\varphi$-absolutely continuous}, namely which are locally (uniformly) $\varphi$-absolutely
continuous in the sense of Tonelli. This means that, for every
$I=\prod_{i=1}^N [a_i,b_i]\subset \Rp^N$ and for every
$j=1,2,\dots,N$, the $j-$th sections of $f$, $f({\tt x}'_j,\cdot):[a_j,b_j]\freccia \erre$, are (uniformly)
$\varphi$-absolutely continuous for almost every ${\tt x}'_j\in
[{\tt a}'_j,{\tt b}'_j]$ (see, e.g., \cite{BAI,DAR}), i.e., there exists $\lambda>0$ such that, for every
$I=\prod_{i=1}^N [a_i,b_i]\subset \Rp^N$, the following property holds:
\vskip0.3cm

\noindent {\it for every
$\varepsilon>0,$ there exists $\delta>0$ such that
$$
 \sum_{i=1}^n \varphi(\lambda|
 f({\tt x}'_j,\beta^i)-f({\tt x}'_j\alpha^i)|)<\varepsilon,
 $$
 for a.e. ${\tt x}'_j\in [{\tt a}'_j,{\tt b}'_j]$ and for all finite collections of non-overlapping intervals
 $[\alpha^i,\beta^i]\subset [a_j,b_j]$, $i=1,\dots,n$, for which
 $$
 \sum_{i=1}^n \varphi(\beta^i-\alpha^i)<\delta.
 $$}

The space $AC^{\varphi}(\Rp^N)$ of the {\it $\varphi-$absolutely continuous} functions will be the space of the functions $f\in {\mathcal M}$ which are of bounded $\varphi-$variation and locally $\varphi-$absolutely continuous on $\Rp^N$.

Strictly related to convergence problems is the notion of modulus of smoothness: in this paper we will use the concept of 
{\it $\varphi-$modulus of smoothness} of $f\in
BV^{\varphi}(\Rp^N)$ defined as
$$
\omega^{\varphi}(f,\delta):= \sup_{|{\bf 1}-{\tt t}|\le\delta} V^{\varphi}[\tau_{\tt t} f-f],
$$
$0<\delta <1$, which is the natural generalization, in the present setting of $BV^{\varphi}(\Rp^N)$, of the classical modulus of continuity (see, e.g., \cite{MU, BMV, AVI5,AVI7}). 
Here $(\tau_{\tt t} f)({\tt s}):= f({\tt st}),$ for every
${\tt s}, {\tt t} \in\Rp^N,$ is the homothetic operator, ${\bf 1}:=(1,\dots,1)$ is the unit vector of $\R^N_+$ and ${\tt s t}:=(s_1 t_1,\dots,s_N t_N)$, ${\tt s,t}\in\R^N_+$. 

\vskip0.3cm

The class of Mellin integral operators that we study is the following:
$$
(T_wf)({\tt s}) = \int_{\Rp^N} K_w ({\tt t}) f({\tt st}) \langle{\tt t}\rangle^{-1} d{\tt t},~~w>0,~ {\tt s} \in \Rp^N, \eqno {\rm(I)}
$$
for $f\in D,$ where $D$ denotes the space of $f:\R^N_+\freccia \R$ for which $(T_w f)({\tt s})$ exists and is finite for every ${\tt s}\in \R^N_+$, $w>0$ (domain of the operators). We remark that $D$ contains a large class of functions, among them, for example, all the bounded functions or, in case of bounded kernels $\{K_w\}_{w>0}$, all the $L^1_{\mu}(\R^N_+)$-functions. \\ Throughout all the paper we will assume that the functions that we consider belong to the domain $D$, so that $(T_w f)({\tt s})$ is well defined for every ${\tt s}\in \R^N_+$, $w>0$. 

As concerns the kernel functions $\{K_w\}_{w>0}$, we assume that:

\begin{description}

\item[${\bf K_w.1)}$] $K_w : \Rp^N\rightarrow \erre$ is a
measurable function such that $K_w \in L^1_{\mu}(\Rp^N),$ $\Vert
K_w\Vert_{L^1_{\mu}}\le A$ for an absolute constant $A>0$ and
$\displaystyle\int_{\Rp^N} K_w ({\tt t}) \langle{\tt t}\rangle^{-1} d{\tt t} =1,$ for every $w>0$;

\item[${\bf K_w.2)}$] for every fixed $0<\delta<1$,
$\displaystyle\int_{|{\bf 1}-{\tt t}|> \delta} |K_w({\tt t})| \langle{\tt t}\rangle^{-1} d{\tt t} \freccia 0$, as $w\freccia +\infty$, 

\end{description}
i.e.,  $\{K_w\}_{w>0}$ is an approximate identity (see, e.g., \cite{BUNE}). We will say that $\{K_w\}_{w>0} \subset \Ks_w$ if $K_w.1)$ and $K_w.2)$ are fulfilled.

\section{Main convergence results}\label{sec-est}

The first result is an estimate for the family of integral operators (I), which shows that our operators map $BV^{\varphi}(\Rp^N)$ into itself.

\begin{prop}\label{prop-est-1}
Let $f\in BV^{\varphi}(\Rp^N)$ and let $\{K_w\}_{w>0}$ be such that $K_w.1)$ holds. Then there exists $\lambda>0$ such that
\begin{equation}\label{varTw-1}
V^{\varphi}[\lambda (T_w f)] \le V^{\varphi}[ \zeta f],
\end{equation}
where $\zeta>0$ is the constant for which $V^{\varphi}[\zeta f] <+\infty$. Therefore, for every $w>0$, $T_w: BV^{\varphi}(\Rp^N) \freccia BV^{\varphi}(\Rp^N)$.
\end{prop}

\noindent{\bf Proof.}\ Let us fix an interval $I=\prod_{i=1}^N[a_i,b_i]\subset \Rp^N$ and a partition of $I$, $\{J_1,\dots,J_m\}$, 
with $J_k=\prod_{j=1}^N[^{(k)}a_j,\ ^{(k)}b_j]$, $k=1,\dots, m$. Let
$\{s_j^o= ~ ^{(k)}a_j,\dots,$ $s_j^{\nu}=~ ^{(k)}b_j\}$ be a partition of the
interval $[^{(k)}a_j,^{(k)}b_j]$, for every $j=1,\dots N$, $k=1,\dots m$. Then, for every $\lambda>0$, ${\tt s}'_j\in I'_j$,
\begin{eqnarray*}
S_j&:=&
\displaystyle\sum_{\mu=1}^{\nu} \varphi(\lambda|(T_w
f)({\tt s}_j',s^{\mu}_j)-(T_w f)({\tt s}_j',s^{\mu-1}_j)|)\\ &=&
\sum_{\mu=1}^{\nu} \varphi\left(\lambda\left|\int_{\Rp^N} K_w({\tt
t}) f({\tt s}'_j{\tt t}'_j,s^{\mu}_jt_j) \langle{\tt t}\rangle^{-1} d{\tt t} +\right.\right. \\ &&  \left.\left.  - \int_{\Rp^N} K_w({\tt
t}) f({\tt s}'_j{\tt t}'_j,s^{\mu-1}_jt_j)  \langle{\tt t}\rangle^{-1} d{\tt t} \right| \right)\\
&\le& \sum_{\mu=1}^{\nu} ~\varphi\left(\lambda\int_{\Rp^N} |K_w({\tt t})||f({\tt s}'_j{\tt t}'_j,s^{\mu}_jt_j) -
f({\tt s}'_j {\tt t}'_j,s^{\mu-1}_j t_j)| \langle{\tt t}\rangle^{-1} d{\tt t} \right).
\end{eqnarray*}
Using Jensen's inequality and assumption $K_w.1)$,
\begin{eqnarray*}
S_j &\le&
A^{-1}~\int_{\Rp^N}| K_w({\tt t})| \sum_{\mu=1}^{\nu}\varphi\Big(\lambda A
|f({\tt s}'_j {\tt t}'_j,s^{\mu}_j t_j) - f({\tt s}'_j {\tt t}'_j,s^{\mu-1}_j t_j)| \Big) \langle{\tt t}\rangle^{-1} d{\tt t} \\
&\le& A^{-1}~\int_{\Rp^N} |K_w({\tt t})|~
V^{\varphi}_{[^{(k)}a_j,^{(k)}b_j]}~[\lambda A f({\tt s}'_j {\tt t}'_j,\cdot ~ t_j)] \langle{\tt t}\rangle^{-1} d{\tt t},
\end{eqnarray*}
and therefore, passing to the supremum over all the partitions of $[^{(k)}a_j,^{(k)}b_j]$,
$$
V^{\varphi}_{[^{(k)}a_j,^{(k)}b_j]}[\lambda(T_w f)({\tt s}'_j,\cdot)] \le A^{-1}~\int_{\Rp^N} |K_w({\tt t})|~
V^{\varphi}_{[^{(k)}a_j,^{(k)}b_j]}~[\lambda A f({\tt s}'_j {\tt t}'_j,\cdot ~ t_j)] \langle{\tt t}\rangle^{-1} d{\tt t}.
$$
Then, by the Fubini-Tonelli theorem,
\begin{eqnarray*}
&&\Phi^{\varphi}_j(\lambda (T_w f), J_k) := \int_{^{(k)}{\tt a}'_j}^{^{(k)}{\tt b}'_j}
V^{\varphi}_{[^{(k)}a_j,^{(k)}b_j]}[\lambda(T_w f)({\tt s}'_j,\cdot)] \langle {\tt s}'_j\rangle^{-1} d{\tt s}'_j \\ &&  \le A^{-1}
\int_{^{(k)}{\tt a}'_j}^{^{(k)}{\tt b}'_j} \left\{ \int_{\Rp^N} |K_w({\tt
t})| V^{\varphi}_{[^{(k)}a_j,^{(k)}b_j]} [\lambda A f({\tt s}'_j {\tt t}'_j,\cdot~ t_j)] \langle{\tt t}\rangle^{-1} d{\tt t} \right\} \langle {\tt s}'_j\rangle^{-1} d{\tt s}'_j \\
&& = A^{-1}\int_{\Rp^N} \left\{ \int_{^{(k)}{\tt a}'_j}^{^{(k)}{\tt b}'_j}
V^{\varphi}_{[^{(k)}a_j,^{(k)}b_j]} [\lambda A f({\tt s}'_j {\tt t}'_j,\cdot~ t_j)] \langle {\tt s}'_j\rangle^{-1} d{\tt s}'_j\right\}
|K_w({\tt t})| \langle{\tt t}\rangle^{-1} d{\tt t} \\ &&= A^{-1} \int_{\Rp^N} \Phi^{\varphi}_j~(\lambda A \tau_{\tt t} f, J_k)|K_w({\tt t})| \langle{\tt t}\rangle^{-1} d{\tt t},
\end{eqnarray*}
for every
$j=1,\dots,N$. Now, applying a Minkowski-type inequality, for every $k=1,\dots,m$ there holds:
\begin{eqnarray*}
\Phi^{\varphi}(\lambda (T_w f), J_k) &:=&\left\{\sum_{j=1}^N [\Phi^{\varphi}_j(\lambda (T_w
f),J_k)]^2\right\}^{1\over 2} \\ &\le& A^{-1} \left\{\sum_{j=1}^N \left(
\int_{\Rp^N} \Phi^{\varphi}_j(\lambda A \tau_{\tt t} f, J_k)|K_w({\tt t})| \langle{\tt t}\rangle^{-1} d{\tt t} \right)^2 \right\}^{1\over 2} \\ &\le& A^{-1} \int_{\Rp^N} \left\{
\sum_{j=1}^N [\Phi^{\varphi}_j(\lambda A \tau_{\tt t} f, J_k)]^2 \right\}^{1\over 2}
|K_w({\tt t})| \langle{\tt t}\rangle^{-1} d{\tt t} \\ &=& A^{-1} \int_{\Rp^N} \Phi^{\varphi}( \lambda A \tau_{\tt t} f, J_k) |K_w({\tt t})| \langle{\tt t}\rangle^{-1} d{\tt t}.
\end{eqnarray*}
Summing over $k=1,\dots,m$ and passing to the supremum over
all the partitions $\{J_1,\dots,J_m\}$ of the interval $I,$ we obtain that
\begin{equation}\label{est0}
V^{\varphi}_I[\lambda (T_w f)]  \le A^{-1} \int_{\Rp^N} V^{\varphi}_I[ \lambda A \tau_{\tt t} f ] |K_w({\tt t})| \langle{\tt t}\rangle^{-1} d{\tt t} ,
\end{equation}
and hence, by the arbitrariness of $I\subset \Rp^N$ and by $K_w.1)$,
$$
V^{\varphi}[\lambda (T_w f)]  \le A^{-1}\Vert K_w\Vert_{L^1_{\mu}} V^{\varphi}[\lambda A f] \le V^{\varphi}[\lambda A f].
$$
Therefore the thesis follows for $0<\lambda \le A^{-1}\zeta$, since $V^{\varphi}[\zeta f]<+\infty$. \hfill\qed

\vskip0.3cm

\noindent{\bf Remark 1.}\ We point out that, in case of $\varphi(u)=u$, $u\in\R^+_0$, and non-negative kernels $\{K_w\}_{w>0}$, then $A=\Vert K_w\Vert_{L^1_{\mu}}=1$, $w>0$, and we can take $\lambda =\zeta=1$. Hence the previous result gives a non-augmenting property of $\varphi-$variation.  

\vskip0.3cm

The following estimate of the error of approximation \hbox{$(T_w f-f)$} with respect to the $\varphi-$variation will be crucial for the main convergence result (Theorem \ref{th-conv}).

\begin{prop} \label{prop-est-2}
Let $f\in BV^{\varphi}(\Rp^N)$ and let $\{K_w\}_{w>0}$ be such that $K_w.1)$ is satisfied. Then for every
$\lambda>0,$ $\delta\in ]0,1[$ and $w>0$,
\begin{eqnarray*}
V^{\varphi}[\lambda(T_w f-f)] \le \omega^{\varphi}(\lambda A f,\delta) +A^{-1} V^{\varphi}[2 \lambda A f] \int_{|{\bf 1}-{\tt t}|>\delta} |K_w({\tt t})| \langle{\tt t}\rangle^{-1}  \,d{\tt t}.
\end{eqnarray*}
\end{prop}

\noindent{\bf Proof.}\ Similarly to Proposition \ref{prop-est-1} (following an analogous reasoning for $(T_w f -f)$, instead of $T_w f$, and recalling that $\displaystyle\int_{\R^N_+} K_w({\tt t}) \langle {\tt t}\rangle^{-1} \,d{\tt t}=1$), it is possible to reach an analogous estimate to (\ref{est0}), i.e., for every $\lambda>0$,
\begin{eqnarray*}
V^{\varphi}_I[\lambda (T_w f-f)]  \le A^{-1} \int_{\Rp^N} V^{\varphi}_I[ \lambda A (\tau_{\tt t}f-f )] |K_w({\tt t})| \langle{\tt t}\rangle^{-1} d{\tt t},
\end{eqnarray*}
and hence, for every $\delta\in]0,1[$, 
\begin{eqnarray*}
V^{\varphi}_I[\lambda (T_w f-f)] & \le A^{-1} \Big(\displaystyle\int_{|{\bf 1}-{\tt t}|\le \delta} + \displaystyle\int_{|{\bf 1}-{\tt t}|>
\delta}\Big) V^{\varphi}_I[ \lambda A (\tau_{\tt t}f-f )] |K_w({\tt t})| \langle{\tt t}\rangle^{-1} d{\tt t}.
\end{eqnarray*}
About the second integral, let us recall that, for every $g,h\in BV^{\varphi}(\R^N_+)$, \hbox{$\lambda>0$,} $V^{\varphi}[\lambda(g+h)] \le 
{1\over 2} \Big(V^{\varphi}[2\lambda g] + V^{\varphi}[2\lambda h]\Big)$ (see property (A) in \cite{AVI7} and also Proposition 1 of \cite{A1}). Therefore
\begin{eqnarray*}
V^{\varphi}_I[\lambda (T_w f-f)] & \le& A^{-1} \int_{|{\bf 1}-{\tt t}|\le \delta}  V^{\varphi}_I[ \lambda A \tau_{\tt t}f-f ] |K_w({\tt t})| \langle{\tt t}\rangle^{-1} d{\tt t}  \\  &+&  {A^{-1}\over 2} \int_{|{\bf 1}-{\tt t}|>
\delta} \Big(V^{\varphi}_I[ 2\lambda A \tau_{\tt t}f] + V^{\varphi}_I[ 2\lambda A f] \Big) |K_w({\tt t})| \langle{\tt t}\rangle^{-1} d{\tt t}.
\end{eqnarray*}
Finally, by the arbitrariness of $I\subset \Rp^N$ and $K_w.1)$, we conclude that
\begin{equation}\label{est1}
\begin{split}
V^{\varphi}[\lambda (T_w f-f)]
&\le
A^{-1} \left\{ \int_{|{\bf 1}-{\tt t}|\le \delta}
|K_w({\tt t})| V^{\varphi}[\lambda A |\tau_{\tt t} f-f|] \langle {\tt t}\rangle^{-1} d{\tt t}\right. \\  &+ \left. V^{\varphi}[2\lambda Af] \int_{|{\bf 1}-{\tt t}|>\delta} |K_w({\tt t})| \langle {\tt t}\rangle^{-1}  d{\tt t}\right\} \\ 
&\le
\omega^{\varphi}(\lambda Af,\delta)+A^{-1} V^{\varphi}[2\lambda Af] \int_{|{\bf 1}-{\tt t}|>\delta} |K_w({\tt t})|  \langle {\tt t}\rangle^{-1}  d{\tt t}. 
\end{split}
\end{equation} \hfill\qed

We are now ready to state the main convergence result.

\begin{thm} \label{th-conv}
Let $f\in AC^{\varphi}(\Rp^N)$ and
$\{K_w\}_{w>0} \subset {\cal K}_w$. Then there exists a constant
$\mu>0$ such that
$$
\lim_{w\to+\infty} V^{\varphi}[\mu(T_w f-f)] =0.
$$
\end{thm}

\noindent{\bf Proof.}\
By Proposition \ref{prop-est-2}, for every $\mu>0$, $\delta\in ]0,1[$, 
$$
V^{\varphi}[\mu(T_w f-f)] \le \omega^{\varphi}(\mu A f,\delta)
+ A^{-1} V^{\varphi}[2\mu A f] \int_{|{\bf 1}-{\tt t}|> \delta}  |K_w({\tt t})| \langle {\tt t}\rangle^{-1}   d{\tt t}.
$$
Now, using Theorem 4.3 of \cite{AVI7}, for every fixed
$\varepsilon>0$ there exist $\overline\lambda>0$ and $0< \overline\delta<1$ such that
$\omega^{\varphi}(\overline\lambda f, \overline\delta)<\varepsilon$ if $|{\bf 1}-{\tt t}|\le \overline\delta$. This imples that $\omega^{\varphi}(\mu A f, \overline\delta)<\varepsilon$ for $0<\mu \le A^{-1}\overline\lambda$. Moreover for every $\delta\in ]0,1[$, by $K_w.2)$, $\displaystyle \int_{|{\bf 1}-{\tt t}|> \delta}  |K_w({\tt t})| \langle {\tt t}\rangle^{-1}  
d{\tt t}<\varepsilon$,  for $w>0$ large enough. Finally, $V^{\varphi}[\zeta f]< +\infty$, for some $\zeta>0$, since $f\in BV^{\varphi}(\Rp^N)$. Therefore, if we consider
$0<\mu\le \displaystyle \min\left\{{\overline\lambda \over A}, {\zeta
\over 2A}\right\},$ then
$$
V^{\varphi}[\mu(T_w f-f)] \le  \omega^{\varphi}(\overline\lambda   f, \overline\delta) 
+  \varepsilon A^{-1} V^{\varphi}[\zeta f] \le  \varepsilon\big(1+A^{-1} V^{\varphi}[\zeta f]\big),
$$
for sufficiently large $w>0$. Hence the theorem is proved, since $\varepsilon>0$ is arbitrary.

\hfill\qed

\noindent{\bf Remark 2.}\ We point out that the assumption that $f\in AC^{\varphi}(\R^N_+)$ in Theorem \ref{th-conv} is essential and cannot be relaxed. For example, the result is no more true, in general, if we just assume that $f\in BV^{\varphi}(\R^N_+)$. Indeed, let us consider, for example in the case $N=1$, the function
$$
f(x)=\begin{cases} 0, \ 0<x<1,\\ 1,\ x\ge 1,
\end{cases}
$$
which is of bounded $\varphi-$variation on $\R_+$, but not $\varphi-$absolutely continuous, and the Mellin Gauss-Weierstrass kernels (see, e.g., \cite{BJ1} and \cite{AVI7} for their multidimensional version) defined as $G_w(t)={w\over \sqrt{\pi}} e^{-{w^2}\log^2 t}$, $t>0$, $w>0$. Then $\{G_w\}_{w>0}$ are approximate identities, i.e., $\{G_w\}_{w>0}\subset \Ks_w$,
$$
(T_w f)(s) =
{1\over \sqrt{\pi}} \int_{{w}\log\left({1\over s}\right)}^{+\infty} e^{-u^2}\,du,\ \ s>0,
$$
and therefore $f\in D$ since $(T_w f)(s)<+\infty$, for every $w>0$. Moreover, for every $\mu>0$,
\begin{align*}
V^{\varphi}[\mu(T_w f-f)] & \ge  V^{\varphi}_{]0,1[}[\mu(T_w f-f)]  = \varphi\left(\mu\left | \lim_{s\to 0^+} (T_wf) (s)-\lim_{s\to 1^-}(T_w f)(s)\right|\right)  \\ & =\varphi\left({\mu\over \sqrt{\pi}}  \int_{0}^{+\infty} e^{-u^2}\,du \right) = \varphi\left({\mu\over 2}\right)>0,
\end{align*}
for every $w>0$, and therefore $V^{\varphi}[\mu(T_w f-f)] \nrightarrow 0$, as $w\freccia +\infty$, for every $\mu >0$.

\section{Order of approximation}\label{sec-order}

In this section we will study the problem of the rate of approximation for the family of operators (I). Before giving the main result, we introduce some definitions.

We say that $\{K_w\}_{w>0}$ is an {\it
$\alpha$-singular kernel}, for $\alpha>0$, if
\begin{eqnarray}\label{ass0}
\int_{|{\bf 1}-{\tt t}|> \delta}  |K_w({\tt t})| \langle {\tt t} \rangle^{-1} d{\tt t}
=O(w^{-\alpha}),~{\rm as}~ w \freccia +\infty,
\end{eqnarray}
for every $\delta\in]0,1[$.

As it is usual in this kind of problems, we have to introduce a
Lipschitz class $V^{\varphi}Lip_{N}(\alpha)$ defined as
\begin{equation*}
V^{\varphi}Lip_{N}(\alpha):=\left\{f\in BV^{\varphi}(\Rp^N)
:~\exists \mu>0\ \hbox{s.t.}\ V^{\varphi}\left[\mu \Delta_{\tt t} f
\right]
 =O(|\log {\tt t}|^{\alpha}),\\ \hbox{as}\ |{\bf 1}-{\tt t}|\freccia 0\right\},
\end{equation*}
where $\Delta_{{\tt t}} f({\tt x}):= (\tau_{{\tt t}} f-f)({\tt x})=
f({\tt x t})-f({\tt x}),$ for ${\tt x,t}\in\Rp^N$, and $\log {\tt t}:=(\log t_1, \dots, \log t_N)$.

\begin{thm}\label{th-order}
Let us assume that
 \hbox{$\{K_w\}_{w>0} \subset  \Ks_w$} is an
$\alpha-$singular kernel and that there exists $0<\widetilde \delta<1$ such that
\begin{equation}\label{ass1}
\int_{|{\bf 1}-{\tt t}|\le \widetilde\delta} |K_w({\tt t})| |\log {\tt t}|^{\alpha} \langle {\tt t}\rangle^{-1}d{\tt t} =O(w^{-\alpha}),\
\ \hbox{as}\ w\freccia +\infty.
\end{equation}
Then if $f\in
V^{\varphi}Lip_{N}(\alpha),$ there exists $\lambda>0$ such that
\begin{equation*}\label{est-order}
V^{\varphi}[\lambda(T_w f-f)]=O(w^{-\alpha}),
\end{equation*}
for sufficiently large $w>0$.\\
\end{thm}

\noindent{\bf Proof.}\  By (\ref{est1}) of Proposition \ref{prop-est-2} we have that, for every $\lambda>0$, $\delta\in ]0,1[$ and $w>0$,
\begin{align*}
V^{\varphi}[\lambda(T_w f-f)] &\le   A^{-1}\left\{
\int_{|{\bf 1}-{\tt t}|\le \delta} V^{\varphi}[\lambda A| \tau_{\tt t} f-f|] |K_w({\tt t})| \langle{\tt t}\rangle^{-1} d{\tt t}  \right.
\\ & \left.+ V^{\varphi}[2 \lambda A
f] \int_{|{\bf 1}-{\tt t}|>\delta} |K_w({\tt t})| \langle{\tt t}\rangle^{-1} d{\tt t}\right\}\\ &:= A^{-1}( J_1+J_2).
\end{align*}
Since $f\in V^{\varphi}Lip_{N}(\alpha),$ there exist $N>0$ and $\bar\delta\in ]0,1[$ such that
$
V^{\varphi}[\lambda A| \tau_{\tt t} f-f|] < V^{\varphi}[\mu \Delta_{\tt t} f] \le N |\log {\tt t}|^{\alpha},
$
if $|{\bf 1}-{\tt t}|\le\bar\delta$ and $0<\lambda 
<\mu A^{-1}$. Now, (\ref{ass1}) ensures that, if $0<\delta\le\min\{\widetilde\delta,\bar\delta\},$ then
$$
J_1\le  N \displaystyle\int_{ |{\bf 1}-{\tt t}|\le
\delta} |K_w({\tt t})| |\log {\tt t}|^{\alpha}
\langle{\tt t}\rangle^{-1} d{\tt t}=O(w^{-\alpha}),
$$
for sufficiently large $w>0$.

Finally, there exist
$\bar\lambda>0$, $M>0$ such that $V^{\varphi}[\bar\lambda f]\le M$, since in particular $f\in BV^{\varphi}(\Rp^N)$. Then, if $0<\lambda <\bar \lambda (2A)^{-1}$, by (\ref{ass0}),
$$
J_2\le M \int_{|{\bf 1}-{\tt t}|>\delta} |K_w({\tt t})| \langle{\tt t}\rangle^{-1}d{\tt t} =O(w^{-\alpha}),
$$
for sufficiently large $w>0$. Hence we conclude that 
$$
V^{\varphi}[\lambda(T_w f-f)]  =O(w^{-\alpha}),
$$
as $w\freccia +\infty$, for $0<\lambda<\displaystyle\min\left\{{\mu\over
A},{\bar \lambda\over 2A}\right\}.$ \hfill\qed

\noindent{\bf Remark 3.}\ We point out that it is possible to obtain a more general version of Theorem \ref{th-order} replacing the functions $|\log {\tt t}|^{\alpha}$ and $w^{-\alpha}$ by $\tau({\tt t})$ and $\xi(w)$, respectively, where $\tau:\Rp^N\freccia \errepz$ is a continuous function at ${\tt t}=1$ and such that $\tau({\tt t})=0$ if and only if ${\tt t}={\tt 1}$, and $\xi:\erre^+_0\freccia \erre^+_0$ is such that $\xi(w)\rightarrow 0$ as $w\rightarrow +\infty$.  
The Lipschitz class has to be now defined as
$$
V^{\varphi}Lip_{N}(\tau):=\left\{f\in BV^{\varphi}(\Rp^N)
:~\exists \mu>0\ \hbox{s.t.}\ V^{\varphi}\left[ \mu \Delta_{\tt t} f
\right]
 =O(\tau({\tt t})),\ \ \hbox{as}\ |{\bf 1}-{\tt t}|\freccia 0\right\},
$$
and (\ref{ass1}) has to be replaced by
$$
\int_{|{1-\tt t}|\le \widetilde\delta} |K_w({\tt t})| \tau({\tt t}) \langle {\tt t} \rangle^{-1}
d{\tt t} =O(\xi(w)),\ \ \hbox{as}\ w\freccia +\infty, \eqno(\ref{ass1}')
$$
for some $\widetilde\delta\in]0,1[$. Finally, $\alpha-$singularity becomes now $\xi-$singularity, i.e., 
\begin{eqnarray*}
\int_{|{\bf 1}-{\tt t}|> \delta}  |K_w({\tt t})| \langle {\tt t} \rangle^{-1}
d{\tt t}
=O(\xi(w)),~{\rm as}~ w \freccia +\infty,
\end{eqnarray*}
for every $\delta\in]0,1[$. Then, similarly to Theorem \ref{th-order} it is possible to prove that 
$$
V^{\varphi}[\lambda (T_w f-f)]=O(\xi(w)),
$$
as $w\freccia +\infty$, for $f\in V^{\varphi}Lip_{N}(\tau)$ and assuming that ($\ref{ass1}'$) holds and that the family $\{K_w\}_{w>0}$ is a $\xi-$singular kernel.

\vskip0.5cm

It is not difficult to find examples of kernel functions which fulfill all the assumptions of Theorem \ref{th-order}. For example, in \cite{AVI9} it is proved that the moment-type kernels defined as
$$
M_w({\tt t}) := w^N \langle {\tt t}\rangle^w \chi_{]0,1[^N}({\tt t}),\ \ {\tt t}\in \Rp^N,\ w>0,
$$
satisfy all the previous assumptions.\\
Moreover, in the classical case, as it is well known, an important class of kernels which satisfy all the assumptions for the rate of approximation is given by the Fej\'er-type kernels with finite absolute moments of order $\alpha$ ($\alpha>0$). The same holds in the present setting, where the Fej\'er-type kernels are kernel functions of the form 
\begin{equation}\label{def_Fejer}
K_w({\tt t}) =w^N K({\tt t}^w),\ \ {\tt t}\in \Rp^N,\ w>0,
\end{equation}
where $K\in L^1_{\mu}(\Rp^N)$ is such that $\int_{\R^N_+} K({\tt t}) \langle {\tt t}\rangle^{-1} \,d {\tt t}=1$, and the absolute moments of order $\alpha$ are defined as
\begin{equation*}\label{moments}
m(K,\alpha):=\int_{\Rp^N} |\log {\tt t}|^{\alpha}| K({\tt t})|\langle{\tt t}\rangle^{-1}  \,d{\tt t}.
\end{equation*}
Indeed in \cite{AVI9} the following Proposition is proved:
\begin{prop} \label{prop_Fejer}
Let $\{K_w\}_{w>0}$ be of the form {\rm(\ref{def_Fejer})} and assume that $m(K,\alpha)<+\infty$. Then 
\begin{description}
\item{(a)} $\displaystyle\int_{|{\mathbf 1}-{\tt t}|>\delta}  |K_w( {\tt t})| \langle{\tt t}\rangle^{-1} \,d{\tt t} =O(w^{-\alpha}),$ as $w\freccia +\infty$, for every $\delta \in ]0,1[$;
\item{(b)} $\displaystyle\int_{|{\mathbf 1}-{\tt t}|\le \delta}  |K_w( {\tt t})| |\log{\tt t}|^{\alpha}\langle{\tt t}\rangle^{-1}  \,d{\tt t} =O(w^{-\alpha}),$ as $w\freccia +\infty$, for every $\delta \in ]0,1[$.
\end{description}
\end{prop} 

\vskip0.3cm 

Finally we point out that there are many examples of Fej\'er-type kernels for which the absolute moments are finite. Among them, there are the Mellin-Gauss-Weierstrass kernels (see \cite{AVI9} and also \cite{AVI6}), defined as
$$
G_w({\tt t}):={w^N\over \pi^{N\over 2}} e^{- w^2 |\log {\tt t}|^2}, \ \ {\tt t} \in \Rp^N,\ w>0;
$$
they are of Fej\'er-type and their absolute moments of order $\alpha$ are finite (\cite{AVI9}). 
Another example are the Mellin Picard kernels, defined as 
\begin{equation*}
P_w({\tt t}):={w^N\over 2\pi^{N\over 2}} {\Gamma({N\over 2})\over \Gamma(N)} e^{-w|\log {\tt t}|}, \ \ t\in \Rp^N,\ w>0,
\end{equation*}
where $\Gamma$ is the Euler function.
Such kernel functions are setted in the frame of $\Rp^N$ from the classical Picard kernels (see, e.g., \cite{BUNE,BABUSTVI,AVI2}), and they are an example of kernels which fulfill all the previous assumptions. First of all they are of Fej\'er-type since
 $P_w({\tt t})=w^N P({\tt t}^w)$ with $P({\tt t})= {\Gamma({N\over 2})\over 2 \pi^{N\over 2}\Gamma(N)}e^{-|\log {\tt t}|}$, ${\tt t}\in\Rp^N$ and $\int_{\Rp^N} P({\tt t})  \langle {\tt t}\rangle^{-1} \,d{\tt t }=1$. Indeed 
\begin{align*} 
I &:=\int_{\Rp^N} P({\tt t})  \langle {\tt t}\rangle^{-1} \,d{\tt t }  = {\Gamma({N\over 2})\over \Gamma(N)} \int_{\Rp^N} e^{-|\log {\tt t}|} \langle {\tt t}\rangle^{-1} \,d {\tt t} = {\Gamma({N\over 2})\over 2\pi^{N\over 2} \Gamma(N)} \int_{\R^N} e^{-|{\tt u}|} \,d {\tt u}. 
\end{align*}
Passing to polar coordinates
$$
\begin{cases}
u_1=\rho \sin \phi_1 \dots \sin\phi_{N-1},
\\ u_2=\rho \sin \phi_1 \dots \cos\phi_{N-1},\\
\dots\\
u_N=\rho\cos \phi_{1},
\end{cases}
$$
and taking into account that, by the Wallis'integrals formula, $\int_0^{\pi\over 2} \sin^n x\,dx={\Gamma\left({n+1\over 2}\right) \Gamma\left({1\over 2}\right) \over 2\Gamma\left({n+2\over 2}\right)}$, then 
\begin{align*}
\int_{\R^N} e^{-|{\tt u}|} \,d{\tt u} & = 
\int_0^{+\infty} e^{-\rho} \rho^{N-1}\,d\rho \int_0^{\pi} \sin^{N-2} \phi_1 \,d\phi_1 \dots \int_0^{2\pi}   \,d\phi_{N-1} \\ & = 2^{N-1}\pi \Gamma(N) \int_0^{\pi\over 2} \sin^{N-2} \phi_1 \,d\phi_1 \dots \int_0^{\pi \over 2}\sin \phi_{N-2} \,d\phi_{N-2} \\ &= 2\pi^{N\over 2} {\Gamma(N) \over \Gamma\left({N\over 2}\right)},
\end{align*} 
and so $I=1$.

Moreover, putting ${\tt u}=\log{\tt t}$, 
\begin{align*}
m(P,\alpha) &= \int_{\Rp^N} |\log{\tt t}|^{\alpha} |P({\tt t})| \langle {\tt t}\rangle^{-1} \,d{\tt t} = 
 {\Gamma({N\over 2})\over 2 \pi^{N\over 2}\Gamma(N)} \int_{\Rp^N} |\log{\tt t}|^{\alpha} e^{-|\log{\tt t}|} \langle {\tt t}\rangle^{-1} \,d{\tt t} \\ &  = {\Gamma({N\over 2})\over 2 \pi^{N\over 2}\Gamma(N)} \int_{\R^N} |{\tt u}|^{\alpha} e^{-|{\tt u}|} \,d{\tt u} <+\infty,
\end{align*}
and hence $\{P_w\}_{w>0}$ are an example of kernel functions to which our results can be applied.

\vskip0.5cm


\end{document}